\documentclass[12pt]{article}
\usepackage{amsmath,amssymb,amsthm}
\usepackage[margin=3.2cm]{geometry}
\usepackage{titlesec,hyperref}
\usepackage{fancyhdr}  

\newtheorem{theo}{Theorem}[section]
\newtheorem{lemm}[theo]{Lemma}

\numberwithin{equation}{section}
\newcommand{\pa}{\partial}
\newcommand{\bbal}{\begin{align*}}
\newcommand{\ri}{\rightarrow}

\begin{document}
\title{Non-uniform dependence for a generalized Degasperis-Procesi equation}
\author{ \footnote{Email Address: guishaohui@126.com (S. Gui),\quad  lijinlu@gnnu.cn (J. Li), \quad  mathzwp2010@163.com (W. Zhu, Corresponding author).}
Shaohui $\mbox{Gui}^{1}$   \quad Jinlu $\mbox{Li}^{1}$  \quad and \quad Weipeng $\mbox{Zhu}^{2}$ \\
   \small $^1\mbox{School}$  of Mathematics and Computer Sciences, \\
\small Gannan Normal University, Ganzhou 341000, China\\
 \small $^2\mbox{Department}$ of Mathematics, Sun Yat-sen University, Guangzhou, 510275, China\\
}
\date{}
\maketitle

\begin{abstract}
In the paper, we consider the Cauchy problem for a generalized Degasperis-Procesi equation. We prove that the data-to-solution map is not uniformly continuous.
\vspace*{1em}

\noindent {\bf Key Words:} A generalized Degasperis-Procesi equation; non-uniform dependence.

\vspace*{1em}

\noindent {\bf Mathematics Subject Classification (2010):} 35Q53
\end{abstract}

\section{Introduction}

In this paper, we study a generalized Degasperis-Procesi equation introduced by Novikov in \cite{n1}:
\begin{equation}\label{E01}
(1-\partial^2_x)u_t=\partial_x(2-\partial_x)(1+\partial_x)u^2.
\end{equation}
It was showed in \cite{n1} that Eqs. (\ref{E01}) possesses a hierarchy of local higher symmetries and the first non-trivial one is $u_{\tau}=\partial_x[(1-\partial_xu)]^{-1}$.

Eqs. (\ref{E01}) belongs to the following class \cite{n1}:
\begin{align}\label{E02}
(1-\partial^2_x)u_t=F(u,u_x,u_{xx},u_{xxx}),
\end{align}
which has attracted much attention on the possible integrable members of (\ref{E02}).

The first well-known integrable member of (\ref{E02}) is the Camassa-Holm (CH) equation \cite{Camassa}:
\begin{align*}
(1-\partial^2_x)u_t=-(3uu_x-2u_{x}u_{xx}-uu_{xxx}).
\end{align*}
The CH equation can be regarded as a shallow water wave equation \cite{Camassa,Camassa.Hyman,Constantin.Lannes}.  It is completely integrable \cite{Camassa,Constantin-P}, has a bi-Hamiltonian structure \cite{Constantin-E,Fokas}, and admits exact peaked solitons of the form $ce^{-|x-ct|}$, $c>0$, which are orbitally stable \cite{Constantin.Strauss}. It is worth mentioning that the peaked solitons present the characteristic for the traveling water waves of greatest height and largest amplitude and arise as solutions to the free-boundary problem for incompressible Euler equations over a flat bed, cf. \cite{Constantin2,Constantin.Escher4,Constantin.Escher5,Toland}. The local well-posedness for the Cauchy problem of the CH equation in Sobolev spaces and Besov spaces was discussed in \cite{Constantin.Escher,Constantin.Escher2,d1,d3,Guillermo,Li-Yin}. It was shown that there exist global strong solutions to the CH equation \cite{Constantin,Constantin.Escher,Constantin.Escher2} and finite time blow-up strong solutions to the CH equation \cite{Constantin,Constantin.Escher,Constantin.Escher2,Constantin.Escher3}. The existence and uniqueness of global weak solutions to the CH equation were proved in \cite{Constantin.Molinet, Xin.Z.P}. The global conservative and dissipative solutions of CH equation were discussed in \cite{Bressan.Constantin,Bressan.Constantin2}.

The second well-known integrable member of (\ref{E02}) is the Degasperis-Procesi (DP) equation \cite{D-P}:
\begin{align*}
(1-\partial^2_x)u_t=-(4uu_x-3u_{x}u_{xx}-uu_{xxx}).
\end{align*}
The DP
equation can be regarded as a model for nonlinear shallow water
dynamics and its asymptotic accuracy is the same as the
CH shallow water equation \cite{D-G-H}. The DP equation is integrable and has a bi-Hamiltonian structure \cite{D-H-H}. An inverse
scattering approach for computing $n$-peakon solutions to the
DP equation was presented in \cite{Lu-S}. Its
traveling wave solutions were investigated in \cite{Le,V-P}. The
local well-posedness of the Cauchy problem of the DP equation in Sobolev spaces and Besov spaces was established in
\cite{G-L,H-H,y1}. Similar to the CH equation, the
DP equation has also global strong solutions
\cite{L-Y1,y2,y4} and finite time blow-up solutions
\cite{E-L-Y1, E-L-Y,L-Y1,L-Y2,y1,y2,y3,y4}. In addition, it has global weak
solutions \cite{C-K,E-L-Y1,y3,y4}. Although the DP equation is similar to the
CH equation in several aspects, these two equations are
truly different. One of the novel features of the DP equation
different from the CH equation is that it has not only
peakon solutions \cite{D-H-H} and periodic peakon solutions
\cite{y3}, but
also shock peakons \cite{Lu} and the periodic shock waves \cite{E-L-Y}.\\

The third well-known integrable member of (\ref{E02}) is the Novikov equation \cite{n1}£º
\begin{align*}
(1-\partial^2_x)u_t=3uu_{x}u_{xx}+u^2u_{xxx}-4u^2u_x.
\end{align*}
The most difference between the Novikov equation and the CH and DP equations is that the former one has cubic nonlinearity and the latter ones have quadratic nonlinearity. It was showed that the Novikov equation is integrable, possesses a bi-Hamiltonian structure, and admits exact peakon solutions $u(t,x)=\pm\sqrt{c}e^{|x-ct|}$, $c>0$ \cite{Hone}. The local well-posedness for the Novikov equation in Sobolev spaces and Besov spaces was studied in \cite{Wu.Yin2,Wu.Yin3,Wei.Yan,Wei.Yan2,Wu.Yin3}. The global existence of strong solutions was established in \cite{Wu.Yin2} under some sign conditions and the blow-up phenomena of the strong solutions was shown in \cite{Wei.Yan2}. The global weak solutions for the Novikov equation were discussed in \cite{Wu.Yin}.

The local well-posedness and global existence of strong solutions for the generalized Degasperis-Procesi were studied in \cite{Li-Yin1}.

To our best knowledge, Eqs. (\ref{E01}) can transform into the following equivalent form:
\begin{equation}\label{1.2}
  \begin{cases}
   u_t-2uu_x=\partial_x(1-\partial^2_x)^{-1}(u^2+(u^2)_x),\ t>0,\ x\in\mathbb{R},\\
   u(0,x)=u_0(x),\ \ x\in\mathbb{R}.
  \end{cases}
\end{equation}

Using this result and the method of approximate solutions, we prove the following nonuniform dependence result.
\begin{theo}\label{th1}
If $s>\frac32$, then the data-to-solution map for the generalized Degasperis-Procesi equation defined by the Cauchy problem \eqref{1.2} is not uniformly continuous from
any bounded subset in $H^s$ into $\mathcal{C}([0, T];H^s)$.
\end{theo}

\noindent\textbf{Notations.} Since all spaces of functions are over $\mathbb{R}$, for simplicity, we drop $\mathbb{R}$ in our notations of function spaces if there is no ambiguity.

\section{Proof of the main theorem}

In this section, we will give the proof of the main theorem. Motivated by \cite{H-H}, we  first construct a sequence approximate solutions. Lately, we will show that the distance between approximate solutions and actual solutions is decaying. Finally, we can conclude that the Cauchy problem \eqref{1.2} is not uniformly continuous.

\begin{lemm}(\cite{B.C.D})\label{le1}
For any $s>0$, there exists a positive constant $c=c(s)$ such that
$$||fg||_{H^s}\leq c(||f||_{H^s}||g||_{L^\infty}+||g||_{H^s}||f||_{L^\infty}).$$
For any $s>\frac12$, there exists a positive constant $c$ such that
$$||f||_{L^\infty}\leq c||f||_{H^s}.$$
\end{lemm}

\begin{lemm}(\cite{B.C.D})\label{le1-1}
Let $s>\frac12$. Assume that $f_0\in H^s$, $F\in L^1_T(H^s)$ and $\pa_xv\in L^1_T(H^{s-1})$. If $f\in \mathcal{C}([0,T];H^s)$ solves the following 1-D linear linear transport equation:
\begin{equation*}
\begin{cases}
\pa_tf+v\pa_xf=F,\\
f(0,x)=f_0,
\end{cases}
\end{equation*}
then there exists a positive constant $C=C(s)$ such that
$$||f||_{H^s}\leq e^{CV(t)}\Big(||f_0||_{H^s}+\int^t_0e^{-CV(\tau)}||F(\tau)||_{H^s}\mathrm{d}\tau\Big),$$
where
\begin{equation*}
V(t)=
\begin{cases}
\int^t_0||\pa_xv||_{H^{s-1}} \ \mathrm{d}\tau, \quad s>\frac32,\\
\int^t_0||\pa_xv||_{H^{s}} \ \mathrm{d}\tau, \quad \frac12<s\leq\frac32.
\end{cases}
\end{equation*}
\end{lemm}

\begin{lemm}(\cite{H-H})\label{le2}
Let $\phi\in S(\mathbb{R})$, $\delta>0$ and $\alpha\in \mathbb{R}$. Then for any $s\geq 0$ we have that
\bbal
&\lim\limits_{n\rightarrow \infty}n^{-\frac12\delta-s}||\phi(\frac{x}{n^{\delta}})\cos(nx-\alpha)||_{H^s}=\frac{1}{\sqrt{2}}||\phi||_{L^2},
\\&\lim\limits_{n\rightarrow \infty}n^{-\frac12\delta-s}||\phi(\frac{x}{n^{\delta}})\sin(nx-\alpha)||_{H^s}=\frac{1}{\sqrt{2}}||\phi||_{L^2}.
\end{align*}
\end{lemm}

\begin{lemm}(\cite{Li-Yin1})\label{le2-2}
Let $s>\frac32$ and $u_0\in H^s$. There exists a positive time $T=T(||u_0||_{H^s})$ such that \eqref{E02} has a solution $u\in \mathcal{C}([0,T];H^s)$. Moreover, there exists a constant $C=C(s)>0$ such that
\bbal
||u||_{L^\infty_T(H^r)}\leq C||u_0||_{H^r}, \quad \forall \ r\geq s.
\end{align*}
\end{lemm}

\textbf{Proof of the theorem:} Let $\omega\in\{0,1\}$ and $\varphi$ be a $C_0(\mathbb{R})$ such that
\begin{equation*}
\phi(x)=
\begin{cases}
1, \quad |x|\leq 1,\\
0, \quad |x|\geq 2.
\end{cases}
\end{equation*}
Let $\phi$ be a $C_0(\mathbb{R})$ such that $\phi(x)\varphi(x)=\varphi(x)$. We introduce the following sequence of high frequency approximate solutions:
\begin{align*}
u^{h,n}_{\omega}=n^{-\frac\delta2-s}\varphi(\frac{x}{n^{\delta}})\cos(nx+2\omega t).
\end{align*}
We also define the solution $u^{\ell,n}_{\omega}$ which satisfies the following equation with the low frequency initial data:
\begin{equation*}\begin{cases}
\pa_tu^{\ell,n}_{\omega}-2u^{\ell,n}_{\omega}\pa_xu^{\ell,n}_{\omega}=\partial_x(1-\partial^2_x)^{-1}[(u^{\ell,n}_{\omega})^2+\pa_x(u^{\ell,n}_{\omega})^2],\\
u^{\ell,n}_{\omega}(0,x)= \omega n^{-1}\phi(\frac{x}{n^\delta}).
\end{cases}\end{equation*}
Since
\bbal
&\pa_tu^{h,n}_{\omega}=-2\omega n^{-\frac\delta2-s}\varphi(\frac{x}{n^{\delta}})\sin(nx+2\omega t)=-2 n u^{\ell,n}_{\omega}(0,x)n^{-\frac\delta2-s}\varphi(\frac{x}{n^{\delta}})\sin(nx+2\omega t),
\\&2u^{\ell,n}_{\omega}\pa_xu^{h,n}_{\omega}=-2nu^{\ell,n}_{\omega} n^{-\frac\delta2-s}\varphi(\frac{x}{n^{\delta}})\sin(nx+2\omega t)+2u^{\ell,n}_{\omega} n^{-\frac{3\delta}{2}-s}\pa_x\varphi(\frac{x}{n^{\delta}})\cos(nx+2\omega t),
\end{align*}
then we can find that
\bbal
\pa_tu^{h,n}_{\omega}-2u^{\ell,n}_{\omega}\pa_xu^{h,n}_{\omega}=&2n [u^{\ell,n}_{\omega}(t,x)-u^{\ell,n}_{\omega}(0,x)]n^{-\frac\delta2-s}\varphi(\frac{x}{n^{\delta}})\sin(nx+2\omega t)\\&-2u^{\ell,n}_{\omega}(t,x) n^{-\frac{3\delta}{2}-s}\pa_x\varphi(\frac{x}{n^{\delta}})\cos(nx+2\omega t).
\end{align*}
Letting $U^n_{\omega}=u^{h,n}_{\omega}+u^{\ell,n}_{\omega}$, we obtain $U^n_{\omega}$ satisfies the following equation:
\bbal
\pa_tU^n_{\omega}-2U^n_{\omega}\pa_xU^n_{\omega}=&-2u^{h,n}_{\omega}\pa_xu^{\ell,n}_{\omega}-2u^{h,n}_{\omega}\pa_xu^{h,n}_{\omega}
\\&+2n [u^{\ell,n}_{\omega}(t,x)-u^{\ell,n}_{\omega}(0,x)]n^{-\frac\delta2-s}\varphi(\frac{x}{n^{\delta}})\sin(nx+2\omega t)\\&-2u^{\ell,n}_{\omega}(t,x) n^{-\frac{3\delta}{2}-s}\pa_x\varphi(\frac{x}{n^{\delta}})\cos(nx+2\omega t)
\\&+\partial_x(1-\partial^2_x)^{-1}[(u^{\ell,n}_{\omega})^2+\pa_x(u^{\ell,n}_{\omega})^2].
\end{align*}
Now, letting $V^n_{\omega}$ be the solution of the Cauchy problem for the equation:
\begin{equation*}\begin{cases}
\pa_tV^n_{\omega}-2V^n_{\omega}\pa_xV^n_{\omega}=\partial_x(1-\partial^2_x)^{-1}[(V^n_{\omega})^2+\pa_x(V^n_{\omega})^2],\\
V^n_{\omega}(0,x)=n^{-\frac\delta2-s}\varphi(\frac{x}{n^{\delta}})\cos(nx)+ \omega n^{-1}\phi(\frac{x}{n^\delta}).
\end{cases}\end{equation*}
Denoting $W^n_{\omega}=U^n_{\omega}-V^n_{\omega}$, it easy to show that
\bbal
&\qquad \pa_tW^n_{\omega}-2U^n_\omega\pa_xW^n_{\omega}-2W^n_\omega \pa_xV^n_\omega\\&=-2u^{h,n}_{\omega}\pa_xu^{\ell,n}_{\omega}-2u^{h,n}_{\omega}\pa_xu^{h,n}_{\omega}
\\&+2n [u^{\ell,n}_{\omega}(t,x)-u^{\ell,n}_{\omega}(0,x)]n^{-\frac\delta2-s}\varphi(\frac{x}{n^{\delta}})\sin(nx+2\omega t)\\&-2u^{\ell,n}_{\omega}(t,x) n^{-\frac{3\delta}{2}-s}\pa_x\varphi(\frac{x}{n^{\delta}})\cos(nx+2\omega t)
\\&+\partial_x(1-\partial^2_x)^{-1}[(u^{\ell,n}_{\omega})^2+\pa_x(u^{\ell,n}_{\omega})^2]
\\&-\partial_x(1-\partial^2_x)^{-1}[(V^n_{\omega})^2+\pa_x(V^n_{\omega})^2]:=\sum^5_{i=1}I_i.
\end{align*}
According to Lemmas \ref{le2}-\ref{le2-2}, we have for $\sigma\geq 0$
\begin{align}\label{estimate1}
||u^{h,n}_{\omega}(t)||_{H^{\sigma}}\leq Cn^{\sigma-s}, \quad ||u^{\ell,n}_{\omega}(t)||_{H^\sigma}\leq C||u^{\ell,n}_{\omega}(t)||_{H^s}\leq C||u^{\ell,n}_{\omega}(0)||_{H^s} \leq Cn^{\frac12\delta-1}.
\end{align}
Choosing $\delta>0$ such that $s-1-\delta>\frac12$, we have $||f||_{L^\infty}\leq ||f||_{H^{s-1-\delta}}$. Therefore, by Lemma \ref{le1} and \eqref{estimate1}, we have
\bbal
||I_1||_{H^{s-1}}&\leq C(||u^{h,n}_{\omega}||_{H^{s-1}}||\pa_xu^{\ell,n}_{\omega}||_{H^{s-1}}+||u^{h,n}_{\omega}||_{H^{s-1-\delta}}||\pa_xu^{h,n}_{\omega}||_{H^{s-1}}
\\&\qquad +||u^{h,n}_{\omega}||_{H^{s-1}}||\pa_xu^{h,n}_{\omega}||_{H^{s-1-\delta}})
\\&\leq Cn^{\frac12\delta-2}+Cn^{-1-\delta},
\end{align*}

\bbal
||I_2||_{H^{s-1}}&\leq Cn\int^t_0||\pa_\tau u^{\ell,n}_{\omega}(\tau)||_{H^{s-1}}\mathrm{d} \tau \cdot ||n^{-\frac\delta2-s}\varphi(\frac{x}{n^{\delta}})\sin(nx-2\omega t)||_{H^{s-1}}
\\&\leq C(||u^{\ell,n}_{\omega}||_{H^s}||\pa_xu^{\ell,n}_{\omega}||_{H^s}+||\partial_x(1-\partial^2_x)^{-1}[(u^{\ell,n}_{\omega})^2+\pa_x(u^{\ell,n}_{\omega})^2]||_{H^{s-1}})
\\&\leq C||u^{\ell,n}_{\omega}||^2_{H^s}+||\pa_xu^{\ell,n}_{\omega}||^2_{H^s}) \leq Cn^{\delta-2},
\end{align*}

\bbal
||I_3||_{H^{s-1}}\leq C||n^{-\frac{3\delta}{2}-s}\pa_x\varphi(\frac{x}{n^{\delta}})\cos(nx-2\omega t)||_{H^{s-1}}\leq Cn^{-1-\delta},
\end{align*}

\bbal
||I_4||_{H^{s-1}}\leq C ||\partial_x(1-\partial^2_x)^{-1}[(u^{\ell,n}_{\omega})^2+\pa_x(u^{\ell,n}_{\omega})^2]||_{H^{s-1}} \leq Cn^{\delta-2}
\end{align*}

\bbal
||I_5||_{H^{s-1}}\leq ||\partial_x(1-\partial^2_x)^{-1}[(V^n_{\omega})^2+\pa_x(V^n_{\omega})^2]||_{H^{s-1}}\leq Cn^{\delta-2}.
\end{align*}
Now, setting $\delta<\min\{s-\frac32,1\}$, it follows from Lemma \ref{le1-1} that
\bbal
||W^n_{\omega}||_{H^{s-1}}\leq Cn^{-1}(n^{-\delta}+n^{\delta-1}).
\end{align*}
According to Lemms \ref{le2-2}, we have
\bbal
||W^n_{\omega}(t)||_{H^{s+1}}&\leq ||U^n_{\omega}(t)||_{H^{s+1}}+||V^n_{\omega}(t)||_{H^{s+1}}\\
&\leq ||u^{h,n}_{\omega}(t)||_{H^{s+1}}+C(||u^{\ell,n}_{\omega}(0)||_{H^{s+1}}+||V^n_{\omega}(0)||_{H^{s+1}})\\
&\leq Cn,
\end{align*}
which implies
\begin{align}\label{estimate2}
||W^n_{\omega}||_{H^{s}}\leq C||W^n_{\omega}||^{\frac12}_{H^{s-1}}||W^n_{\omega}||^{\frac12}_{H^{s+1}}\leq C(n^{-\delta}+n^{\delta-1})^{\frac12}.
\end{align}
Then, combining \eqref{estimate1} and \eqref{estimate2}, we have
\begin{align}\label{estimate3}\begin{split}
||V^n_1(t)-V^n_0(t)||_{H^s}&\geq ||U^n_1(t)-U^n_0(t)||_{H^s}-C\varepsilon_n
\\&\geq ||u^{h,n}_1(t)-u^{h,n}_0(t)||_{H^s}-C\varepsilon'_n
\\&\geq 2|\sin t|\cdot||n^{-\frac\delta2-s}\varphi(\frac{x}{n^{\delta}})\sin(nx+t)||_{H^s}-C\varepsilon'_n,
\end{split}\end{align}
where
\bbal
\varepsilon_n=(n^{-\delta}+n^{\delta-1})^{\frac12}, \qquad \varepsilon'_n=(n^{-\delta}+n^{\delta-1})^{\frac12}+n^{\frac12\delta-1}.
\end{align*}
Letting $n$ go to $\infty$ and using Lemma \ref{le2}, \eqref{estimate2}-\eqref{estimate3}, it follows that
\begin{align}\label{estimate4}
\liminf_{n\ri \infty}||V^n_1(t)-V^n_0(t)||_{H^s}&\geq \sqrt{2}||\varphi||_{L^2}|\sin t|.
\end{align}
Noticing that $|\sin t|=\sin t$ when $t\in[0,\pi]$, then $\sin t>0$ in an interval $(0,t_0)$ for some $0<t_0<\pi$. This together with the fact that
\bbal
\lim_{n\ri \infty}||V^n_1(0)-V^n_0(0)||_{H^s}\leq ||n^{-1}\phi(\frac{x}{n^\delta})||_{H^s}\leq n^{\delta-1}\rightarrow 0, \quad \mathrm{as} \ n \rightarrow \infty,
\end{align*}
complete the proof of Theorem \ref{th1}. \\

\noindent\textbf{Acknowledgements.} This work was partially supported by NSFC (No.11671407).

\phantomsection

\end{document}